\documentstyle[12pt]{article}   
\usepackage{latexsym,amsfonts}   
\title{On Clifford representation of Hopf algebras and Fierz identities.} 
\author{       
Suemi Rodr\'\i guez-Romo\\Centre of Theoretical Research\\
National University of M\'exico, Campus Cuautitl\'an\\
Apdo. Postal 95, Unidad Militar, Cuautitl\'an Izcalli\\
Estado de M\'exico, 54768 M\'exico.$^{\ast}$}         
\date{}  
\begin{document}   
\maketitle
\renewcommand{\thefootnote}{\fnsymbol{footnote}}
\setcounter{footnote}{-1}
\footnote{$\hspace*{-6mm}^{\ast}$
e-mail: suemi@servidor.unam.mx \\
$\hspace*{1.4cm}$ suemi@fis.cinvestav.mx}
\renewcommand{\thefootnote}{\arabic{footnote}}
\baselineskip1.2cm
Abstract.  We present a short review of the action and coaction 
of Hopf algebras on Clifford algebras as an introduction to physically
meaningful examples. Some $q$-deformed Clifford algebras are studied from this
context and conclusions are derived.\


Key words; Clifford algebra, quantum groups, Fierz identities, multivectorial
Cartan map.

\newpage

\baselineskip1.2cm  

\section{INTRODUCTION.}
\newcommand{\BC}{{\hbox{\rm\vphantom{X}%
                       \hskip 0.25em%
                       \vrule width 0.7pt%
                       \hskip -0.35em C}}}

\newcommand{\cl}{C \kern -0.1em \ell} 
Like a group, a Hopf algebra can act as the algebra of symmetries of another 
algebra (cf. the group of symmetries of a quantum system). In this paper we
explore the possible quantum Clifford action (coaction) on Clifford algebras. 
Hopf algebras can act on other structures in a variety of ways; i.e., $H$ can 
act as an algebra on a vector space. Besides, we can think of quantization as 
a procedure that replaces the classical algebra of observables by a 
non-commutative quantum algebra of observables \cite{Dr}. It is expected that 
even using non-commutative geometry, one might nevertheless extend our regular 
notions of symmetry to the quantum world. The simplest non-commutative 
geometries that have been studied are Clifford algebras and Hopf algebras 
corresponding to both quantum symmetry and curvature. Thus, it 
is a natural question to inquire about the way Hopf algebras act (and coact) 
on Clifford algebras to understand how quantum structures act (and coact) on 
classical ones as physically meaningful symmetries. Other attempts to link
Clifford algebra to $q$-deformed geometry have been reported recently. See
for example \cite{So}.\

Quantum algebras have already proved useful in the study of quantum spin 
chains and conformal field theories, and may have other important physical 
applications as well.\

\section{Clifford category of representations of a triangular Hopf algebra.}

Let us begin this Section with some formal definitions. A {\it Hopf algebra}
\newline
($H$, $+$, $\cdot$ , $\eta$, $\Delta$, $\epsilon$, $S$: $k$) over $k$, a 
field, is a bialgebra over $k$ equipped with a linear antipode map $S$: 
$H\rightarrow H$ obeying $\cdot (S\otimes id)\circ \Delta$=
$\cdot(id\otimes S)\circ \Delta$=$\eta\circ \epsilon$. Here $\epsilon$ is the 
counit map $\epsilon:H\rightarrow k$, $\Delta $ is the coproduct map 
$\Delta: H\rightarrow H\otimes H$ and $\eta$ is the unit map 
$\eta: k\rightarrow H$. A {\it bialgebra} ($B$, $+$ , $\cdot$, $\eta$, 
$\Delta$, $\epsilon$; $k$) over $k$ is a vector space ($B$, $+$; $k$) over 
$k$ which is both an algebra and a coalgebra, in a compatible way; i.e. 
$\Delta(ab)$=$\Delta(a)\Delta(b)$, $\Delta(1)$=$1\otimes 1$, $\epsilon(ab)$=
$\epsilon(a)$$\epsilon(b)$ and $\epsilon(1)$=$1$. Finally, a {\it coalgebra} 
($\AE$, $+$, $\Delta$, $\epsilon$; $k$) over $k$ is a vector space \newline 
($\AE$, $+$; $k$) over $k$ and a linear map $\Delta$ which is a coassociative  
and for which there exists a linear counit map $\epsilon $. In this paper we
are specially interested in non-commutative and non-cocommutative Hopf 
algebras. This means that neither the product nor the coproduct are 
commutative.\

A {\it representation} of an algebra $H$ is a pair $(\alpha,V)$ where $V$ is 
a vector space and $\alpha$ a linear map $H\otimes V\rightarrow V$, say 
$\alpha(h\otimes v)$=$\alpha_h(v)$ for $h\in H$ and $v\in V$, where
$\alpha_h(\alpha_g(v))$=$\alpha_{hg}(v)$. A left $H$-module is nothing 
other than the vector space on which the algebra $H$ is represented.\

If $\alpha$ is an antirepresentation, namely $\alpha_h(\alpha_g(v))$=
$\alpha_{gh}(v)$, then $\alpha$ is the right action of ${\it H}$ on $V$. In 
other words, $V$ is a right ${\it H}$-module.\ 

The collection of left $H$-modules is denoted by $_{\it H}{\it M}$ and
the collection of right $H$-modules is denoted by ${\it M}_{\it H}$.\

If $H$ is a Hopf algebra then it can act on algebras in such a way as to 
respect the algebra structure. Thus, $H$ acts on $A$ as an algebra ( or $A$ 
is an {\it $H$-module algebra}) if $A$ is an $H$-module and in addition
$$
\alpha(h\otimes ab)=\alpha(h_{(1)}\otimes a)\alpha(h_{(2)}\otimes b),
\;\;\;\alpha(h\otimes 1_{A})=\epsilon1_{A}
$$
for $a,b\in A$ and $h\in H$. Here  ${\bf 1}_A$ is the unit in $A$ and
$\Delta(h)$=$h_{(1)}\otimes h_{(2)}$.\

As an example of a Hopf algebra acting on a Clifford algebra, let us 
consider, formally, $GL_q(2,\BC)$ acting on $\cl_{3,1}$. Here 
$\gamma_{\mu}\in \cl_{3,1}$ such that 
$\left\{\gamma_{\mu},\gamma_{\nu}\right\}$=$2g_{\mu\nu}{\bf 1}$;
$\mu,\nu=0,...,3$ and $g_{\mu\nu}$=diag(-1,1,1,1). On the ohter hand we 
define $GL_q(2,\BC)$ as follows. Let
$a_{ij}\in GL_q(2,\BC)$; $i,j=1,2$, then $a_{11}a_{12}=qa_{12}a_{11}$; 
$a_{11}a_{21}=qa_{21}a_{11}$; $a_{12}a_{21}=a_{21}a_{12}$; 
$a_{12}a_{22}=qa_{22}a_{12}$; $a_{21}a_{22}=qa_{22}a_{21}$ and 
$a_{11}a_{22}-a_{22}a_{11}$=$(q-q^{-1})a_{12}a_{21}$. Besides 
$\Delta(a_{ij})$=$\sum^2_{k=1}a_{ik}\otimes a_{kj}$, $\epsilon(a_{11})$=
$\epsilon(a_{22})$=$1$ and $\epsilon(a_{12})$=$\epsilon(a_{21})$=$0$. Here 
$q$ is any complex number. We adopt the concept of quantum group in a purely
geometrical point of view \cite{Ma}.\

We can see that $q$ plays the rule of a deformation parameter which, in 
the limit $q\rightarrow \pm 1$, includes bose and fermi statistics in Fock 
space. The generators $a_{ij}$ have been identified with creation and 
annihilation operators in $q$-deformed oscillators for a wide variety of
authors. See for example \cite{Ch}. In this case the parameter $q$ allows us 
to generalize statistics to include parfermi, infinite and other cases 
\cite{Sm}. Some applications related to the spectra of triatomic molecules 
and superdeformed nuclei have also been studied.\

Let
$$
\alpha(a_{ij}\otimes \gamma_{\mu}\gamma_{\nu})=
\alpha(a_{ij(1)}\otimes \gamma_{\mu})\alpha(a_{ij(2)}\otimes \gamma_{\nu}) 
\;\;\forall {\mu},{\nu}\;\;\mbox{  and}
$$
$$
\alpha(a_{ij}\otimes \gamma_{\mu}\gamma_{\nu})=\alpha(a_{ij(1)}\otimes 
\gamma_{\mu})\cdot \alpha(a_{ij(2)}\otimes \gamma_{\nu})+
\alpha(a_{ij(1)}\otimes \gamma_{\mu})\wedge \alpha(a_{ij(2)}\otimes 
\gamma_{\nu})=
$$
$$
\alpha(a_{ij(k)}\otimes g_{\mu\nu})+\alpha(a_{ij(k)}\otimes \gamma_{\mu\nu})
$$
for any $\gamma_{\mu},\gamma_{\mu\nu}\in \cl_{3,1}$ and $k=1,2$. From
this, the action of the map on any $\Gamma\in\cl_{3,1}$ follows. A 
concrete map remains to be found.\

In this paper we present two physically meaningful concrete examples related
with this matter.\\

1. {\bf $\cl_{3,1}$ as a ${\it C}_q(3,1)$-module}. In this case we are 
thinking of ${\it C}_q(3,1)$, a $q$-deformed Clifford algebra, as symmetry 
of $\cl_{3,1}$. Following Manin \cite{Ma}, ${\it C}_q(3,1)$ can be defined
by considering a space, the coordinates of which do not commute, as follows. 
Let $\left\{\gamma^{\mu}_q\right\}$
$\in {\it C}_q(3,1)$ then $\gamma^{\mu}_q\gamma^{\nu}_q+q{\hat{\bf R}}^
{\mu\nu}_{\nu'\mu'}\gamma^{\nu'}_q\gamma^{\mu'}_q$=
$q^{-1}QC^{-1 \mu\nu}$, where
$$
\gamma^0_q=
\left(
\begin{array}{cccc}
0  & 0  & q^2 & 0   \\
0  & 0  & 0   & -1  \\
-1 & 0  & 0   & 0   \\
0  & -1 & 0   & 0
\end{array}
\right)
\hspace*{1cm}
\gamma^+_q=
\sqrt{qQ}
\left(
\begin{array}{cccc}
0  & 0   & 0   & 1   \\
0  & 0   & 0   & 0  \\
0  & -1  & 0   & 0   \\
0  & 0  & 0   & 0
\end{array}
\right)
$$
$$
\gamma^-_q=\sqrt{Q}
\left(
\begin{array}{cccc}
0         & 0  & 0   & q^{-3/2}   \\
0         & 0  & 0   & 0          \\
0         & 0  & 0   & 0          \\
-q^{3/2}  & 0  & 0   & 0
\end{array}
\right)
\hspace*{1cm}
\gamma^3_q=
\left(
\begin{array}{cccc}
0   & 0    & q^{-1}+q-q^2   & 0        \\
0   & 0    & 0              & -q^{-2}  \\
-1  & 0    & 0              & 0   \\
0   & q^2  & 0              & 0
\end{array}
\right),
$$
Here $Q=q+q^{-1}$, being q a real number, ${\hat{\bf R}}^{\mu\nu}_{\nu'\mu'}$ 
is the $SO_q(3,1)$-R matrix and $C$ is the following metric \cite {Ca}.
$$
\left(
\begin{array}{cccc}
0   &     0       &     0     &   q^{-1} \\
0   & -1+q^{-2}   & -q^{-1}   &     0    \\
0   &  -q^{-1}    &     0     &     0    \\
q   &      0      &     0     &     0
\end{array}
\right).
$$ 
We define $\alpha$:
$\sum_{\rho}\sum_{\mu}\left(\gamma^{\mu}_q\right)
^{\rho}_{\nu}\otimes\gamma_{\mu}\rightarrow \gamma_{\nu}$. Before $\alpha$ is 
applied $\left\{\gamma_{\mu},\gamma_{\nu}\right\}$=$2g_{\mu\nu}{\bf 1}$ with
$\mu, \nu, \rho=0,+,-,3$ and $g=[g_{\mu\nu}]$=diag(-1,1,1,1). Here the 
indices $0,+,-,3$ correspond to the indices $\mu$ in the
$\gamma^{\mu}_q$ matrices above defined. After applying the map, 
$\left\{\alpha_{\gamma^{\mu}_q},\alpha_{\gamma^{\nu}_q}\right\}$=
$2\alpha_{g^{\mu\nu}}{\bf 1}$, where
$$
\alpha_{g^{\mu\nu}}=
$$
$$
\left(
\begin{array}{cccc}
Qq^{-3}    &   1-q^2       &   q-q^{-1}    &   1-Qq^{-3}+q^{-2}     \\
1-q^2    & Qq+q^4-1      & q-q^2+q^3-q^4 &          qQ              \\
q-q^{-1} & q-q^2+q^3-q^4 &       Q(Q-2q^2) &    1-q^2-q^{-1}-q^{-3} \\
1-Qq^{-3}+q^{-2} &   qQ    &  1-q^2-q^{-1}-q^{-3} & -1+qQ+q^{-3}+q^{-4}
\end{array}
\right).
$$\

From this, the  ${\it C}_q(3,1)$-invariance of the
general Fierz identity \cite{Su1} follows, which comes from the 
multivectorial generalization of the Cartan map \cite{Su3} used to obtain a 
multivectorial Dirac equation \cite{Su4}. In this multivectorial Dirac 
equation, by means of the comodule map $\beta$ defined from the 
specific $\alpha$ in this case (reversing arrows), we can show the
${\it C}_q(3,1)$ symmetry.\

2. {\bf $su(2)$ as a $CH_q(2)$-module algebra}. Let us consider in this case the 
affinization of the Hopf algebra $CH_q(2)$ acting on the $su(2)$ algebra. 
$CH_q(2)$ is the quantum deformation of the Hopf algebra $CH(2)$ \cite{Cu}. 
Here $CH(2)$ is generated by $\Gamma_1$, $\Gamma_2$, $\Gamma_3$ and 
$E_1$, $E_2$, $E_3$ such that $\Gamma^2_{\mu}$=$E_{\mu}$; $\Gamma^2_3=1$; 
$\{\Gamma_{\mu},\Gamma_{\nu}\}=0$, $\mu\neq\nu$; 
$\{\Gamma_{\mu},\Gamma_{3}\}=0$; $[E_{\mu},\Gamma_{\nu}]$=
$[E_{\mu},\Gamma_3]$=$[E_{\mu},E_{\nu}]$=0, $\forall \mu,\nu$;\newline 
$\Delta(E_{\mu})$=$E_{\mu}\otimes {\bf 1}+{\bf 1}\otimes E_{\mu}$;
$S(E_{\mu})=-E_{\mu}$, $\epsilon (E_{\mu})=0$, 
$\Delta(\Gamma_{\mu})$=
$\Gamma_{\mu}\otimes{\bf 1}+\Gamma_3\otimes\Gamma_{\mu}$; $S(\Gamma_{\mu})$=
$\Gamma_{\mu}\Gamma_3$; $\epsilon(\Gamma_{\mu})$=0; $\Delta(\Gamma_3)$=
$\Gamma_3\otimes\Gamma_3$; $S(\Gamma_3)$=$\Gamma_3$; $\epsilon(\Gamma_3)$=1;
where $\mu,\nu$=1,2.\

The quantum deformation of $CH(2)$; i.e. $CH_q(2)$, is carried out by only 
transforming within this, the following terms
$$
\Gamma^2_{\mu}=[E_{\mu}]_q=\frac{q^{E_{\mu}}-q^{-E_{\mu}}}{q-q^{-1}};
\;\;\;
\Delta \Gamma_{\mu}=\Gamma_{\mu}\otimes q^{-E_{\mu}/2} + 
q^{E_{\mu}/2}\Gamma_3\otimes \Gamma_{\mu}.
$$
$\widehat{CH_q(2)}$, the affinization of $CH_q(2)$, is generated by 
$E_{\mu}^{(i)}$, $\Gamma^{(i)}_{\mu}$ $(i=0,1;\mu=1,2)$ and $\Gamma_3$ 
satisfying $CH_q(2)$ for each value of i. Here $q$ is any complex number.\

A two dimensional irrep of $\widehat{CH_q(2)}$ is labelled by 
$(z,\lambda_x,\lambda_y)$$\in {\bf C}^3$ and reads \cite{Cu}; 
$$
\Gamma^{(0)}_x=
\left(\frac{ \lambda^{-1}_x-\lambda_x}{q-q^{-1}}\right)^{1/2}
\left(
\begin{array}{cc}
0 & z^{-1}\\
z & 0
\end{array}
\right)
\;\;,\;\;
\Gamma^{(0)}_y=
\left(\frac{ \lambda^{-1}_y-\lambda_y}{q-q^{-1}}\right)^{1/2}
\left(
\begin{array}{cc}
0  & -iz^{-1}\\
iz &  0
\end{array}
\right),
$$\vspace{.5cm}
$$
\Gamma^{(1)}_x=
\left(\frac{\lambda_x-\lambda^{-1}_x}{q-q^{-1}}\right)^{1/2}
\left(
\begin{array}{cc}
0      & z\\
z^{-1} & 0
\end{array}
\right)
\;\;,\;\;
\Gamma^{(1)}_y=
\left(\frac{\lambda_y-\lambda^{-1}_y}{q-q^{-1}}\right)^{1/2}
\left(
\begin{array}{cc}
0       & -iz\\
iz^{-1} &  0
\end{array}
\right),
$$\vspace{.5cm}
$$
\Gamma_3=
\left(
\begin{array}{cc}
1 & 0\\
0 & -1
\end{array}
\right)
\;,\;q^{E^{(0)}_x}=\lambda^{-1}_x\;,\;q^{E^{(0)}_y}=\lambda^{-1}_y
\;,\;q^{E^{(1)}_x}=\lambda_x\;\;\mbox{  and  }\;\;q^{E^{(1)}_y}=\lambda_y.
$$
Let us define the map $\alpha:\sum_{\rho}\sum_{\mu}
(\Gamma^{(i)\mu})^{\rho}_{\nu}\otimes \sigma_{\mu}\rightarrow \sigma_{\nu}$
where $\mu,\nu, \rho$=$1,2$ and $\alpha: 1\otimes \sigma_3\rightarrow \sigma_3$, 
being $\{\sigma_k\}$ a basis for the Lie algebra of $su(2)$; i.e. the Pauli
matrices. Before this map is applied, $\sigma^2_i=1$ and 
$\{\sigma_i,\sigma_j\}$=$2\delta_{ij}$; $i,j=1,2,3$. After applying the map 
$(\alpha_{\sigma_i})^2$=$0$, $(\alpha_{\sigma_3})^2$=$1$, 
$\{\alpha_{\sigma_1},\alpha_{\sigma_3}\}$=$0$ and 
$\{\alpha_{\sigma_i},\alpha_{\sigma_j}\}$=$4(1\otimes 1)$; $i,j=1,2$.\

From this, it is straightforward to verify the $\widehat{CH}_q(2)$-invariance
of the Clifford product in the extended Cartan map structure \cite {Re1} which 
form a quaternion algebra.\

On the other hand, $H$ acts adjointly on $A$ if $\alpha(h\otimes\Gamma)$=
$h_{(1)}\Gamma Sh_{(2)}$ for any $\Gamma \in A$. Actually in general 
$f(\alpha(h\otimes \Gamma))$=$f(h_{(1)})\Gamma f(Sh_{(2)})$ for any $f$ 
algebra map. Here $S: H\rightarrow H$ is the antipode map. For 
$A=\cl_{3,1}$ and $H=SL_q(2)$ the adjoint action looks like 
$\alpha(a_{ij}\otimes \Gamma)$=$a_{ij(1)}\Gamma S(a_{ij(2)})$, for any
$\Gamma \in \cl_{3,1}$.\

$\cl_Q$ is a $Z_2$-graded algebra; i.e. $\cl_Q$=
$\cl_Q^0\oplus \cl_Q^1$, where $\Lambda^{(r)}(Q)\Lambda^{(s)}(Q)$\newline
=$\Lambda^{(r+s)}(Q)$ and $1_{\cl_Q}$=$\Lambda^{(0)}(Q)$=$k$. Therefore,
any Clifford algebra can be considered as a twisted tensor product of the
$Z_2$-graded algebras. Since $Z_2$ is a unital semigroup then 
$\cl_Q$ does not induce a $k(Z_2)$-module algebra, 
which is the set of functions on $Z_2$ with values in $k$. This is 
because the antipode map cannot be defined.\

Similarly for $A$ being a coalgebra we can require the action of $H$ to respect 
the coalgebra structure; i.e.
$$
\alpha(h\otimes a)_{(1)}\otimes\alpha(h\otimes a)_{(2)}=
\alpha(h_{(1)}\otimes a_{(1)})\otimes\alpha(h_{(2)} \otimes a_{(2)})
\;\;\mbox{ and }
$$
$$
\epsilon\alpha(h\otimes a)=\epsilon(h)\epsilon(a)
$$
Since $\cl_Q$ is not a coalgebra we do not require this 
condition to hold.\

Dual to the notion of modules is the notion of comodules. Since they
are dual their diagrams are obtained by reversing arrows. In other words $V$ 
is a left-${\it H}$-comodule if there is a map 
$\beta:V\rightarrow {\it H}$$\otimes V$ such that 
$(id\otimes \beta)\circ \beta$=\newline$(\Delta \otimes id)\circ\beta$ and
$(\epsilon\otimes id)\circ \beta$$=1_{\it H}\otimes id$. See Figure 1.\

The collection of left $H$-comodules is denoted by $^H{\it M}$ and the 
collection of right comodules is denoted by ${\it M}^{H}$. The three cases of 
quantum symmetry above presented can easily be reversed to obtain comodule 
structures.\

Let ${\it H}$ be a finite dimensional Hopf algebra and ${\it H}^*$ its dual. 
It is well known that there is a one-to-one correspondence between left 
${\it H}$-modules and right ${\it H}^*$-comodules, left ${\it H}$-module 
coalgebras and right ${\it H}^*$-comodule coalgebras, left ${\it H}^*$-module
algebras and right ${\it H}^*$-comodule algebras.\ 

Let $(\alpha_1,V_1)$ and $(\alpha_2,V_2)$ be any two representations of the 
algebra part of $H$ in the vector spaces $V_1$ and $V_2$.\

A {\it category} ${\it C}$ is a collection of objects $Ob({\it C})$, and a 
set $Mor(X,Y)$ for each $X,Y\in Ob({\it C})$. The latter are called morphisms 
and $Mor(X_1,Y_1)$, $Mor(X_2,Y_2)$ are disjoint unless $X_1=X_2$ and 
$Y_1=Y_2$. They should have properties analogous to those of maps from $X$ to 
$Y$ that respect the structure on $X$ and $Y$. The representations of an 
algebra form a category; i.e., the category of the Clifford representations of 
${\it C}_q(3,1)$. Indeed, an object in this category is a pair $(\alpha, V)$. The 
morphisms $Mor((\alpha_1, V_1), (\alpha_2, V_2))$ in this category are the 
interwiners, namely $\phi \in Lin_k(V_1,V_2)$, these are the k-linear maps 
from $V_1$, $V_2$ to k such that
$$
\alpha_2(\phi(v))=\phi(\alpha_1(v))\;,\;\;
\forall v\in V_1
$$
The so called reflection equations for quantum groups can be written in terms 
of interwiners \cite{Az}.\

A map between categories $F:{\it C}_1\rightarrow {\it C}_2$ is called a 
{\it covariant functor} if to each object $X\in Ob({\it C}_1)$ it assigns an object 
$F(X)\in Ob({\it C}_2)$ and to each morphism $\phi\in Mor(X,Y)$ it assigns a 
morphism $F(\phi)\in Mor (F(X),F(Y))$ such that 
$F(\phi_1\circ\phi_2)=F(\phi_1)\circ F(\phi_2)$.
For a {\it contravariant functor}
 \newline $F:{\it C}_1\rightarrow {\it C}_2$, also assigns
an object $F(X)\in Ob({\it C}_2)$ to each object $X\in Ob({\it C}_1)$ but 
assigns an element $F(\phi)\in Mor(F(Y), F(X))$ for each $\phi \in Mor(X,Y)$ 
such that $F(\phi_1\circ\phi_2)=F(\phi_2)\circ F(\phi_1)$.\

On the other hand a {\it natural transformation} $\Phi:F_1\rightarrow F_2$ between
two functors $F_1, F_2:{\it C}_1\rightarrow {\it C}_2$ is a map that assigns 
to each object $X\in Ob({\it C})$ a morphism $\Phi_X\in Mor(F_1(X), F_2(X))$
such that for any morphism \newline $\Phi\in Mor(X,Y)$ in ${\it C}_1$,
$\Phi_Y\circ F_1(\Phi)$=$F_2(\Phi)\circ\Phi_X$
and similarly if $F_1$ and $F_2$ are contravariant.\

A natural transformation $\Phi$ is called a {\it natural equivalence} of 
functors if each map $\Phi_X$ is an isomorphism. The maps $\Phi_X$ in this 
case are also said to be functorial isomorphisms.\

A category $({\it C},\tilde{\otimes}, \underline{1})$, is called monoidal if it
has a product functor $\tilde{\otimes}:{\it C}\times {\it C}$
$\rightarrow {\it C}$ and a unit object 
$\underline{1}\in Ob({\it C})$ such that; \newline
a) The two functors ${\it C}\times {\it C}\times {\it C}\rightarrow {\it C}$
given by $\tilde{\otimes}(\tilde{\otimes})$ and 
$(\tilde{\otimes})\tilde{\otimes}$ are naturally equivalent; i.e. in addition 
to $\tilde{\otimes}$, there are functorial isomorphisms
$$
\Phi_{X,Y,Z}:X\tilde{\otimes}(Y\tilde{\otimes}Z)\rightarrow(X\tilde{\otimes}
Y)\tilde{\otimes}Z. 
$$
b) $X\tilde{\otimes}Y\tilde{\otimes}Z\tilde{\otimes}W$ is associative.\newline
c) Finally, the functors ${\it C}\rightarrow {\it C}$ given by $X\rightarrow
X\tilde{\otimes}\underline{1}$ and $X\rightarrow\underline{1}\tilde{\otimes}X$ 
should be naturally equivalent to the identity functor; i.e. there are 
functorial isomorphisms $X\tilde{\otimes}\underline{1}\rightarrow X$ and 
$\underline{1}\tilde{\otimes}X\rightarrow X$.\
                                                                   
Consider $( _{{\it H}}{\it M}, \tilde{\otimes}, \phi, \underline{1})$,
where $_{{\it H}}{\it M}$ is the category of algebra representations of Hopf 
algebras, $\tilde{\otimes}$ is defined on $_{{\it H}}{\it M}$ by 
$ V_1\tilde{\otimes}V_2=V_1\otimes V_2$ as vector spaces, 
$\alpha(h\otimes(v_1\tilde{\otimes}v_2))$=
$\alpha(h_{(1)}\otimes v_1)\tilde{\otimes}\alpha(h_{(2)}\otimes v_{2})$ 
$\forall$ $v_1\tilde{\otimes}v_2\in V_1\tilde{\otimes}V_2$ and where the unit 
is defined by the trivial representation on $k$. Then $_{{\it H}}{\it M}$ is a 
monoidal category.\

In this case, there is a functor $F:$$ _{H}{\it M}\rightarrow Vec$ that 
assigns to $(\rho,V)$ the vector space $V$ and to a morphism 
$\phi\in Mor((\rho_1 V_1),(\rho_2,V_2))$ the interwiner viewed just as a 
linear map $\phi\in Mor_{Vec}(V_1,V_2)$=$Lin_k(V_1,V_2)$. This is called the 
forgetful functor and respects the monoidal structure.\

We construct the Clifford algebra $\cl_Q$ as the quotient algebra of
$\otimes V$ with respect to the two-sided ideal $I(Q)$ generated by 
the elements \newline $X\otimes X-Q(X)$ where $X\in V$ and $Q$ is a quadratic form 
on $V$ and we consider monoidal tensor categories.\

A tensor category $({\it C}, \tilde{\otimes}, \phi, \Psi, \underline{1})$ is 
a monoidal category $({\it C},\tilde{\otimes},\phi,\underline{1})$ such that 
the two functors ${\it C}\times {\it C}\rightarrow {\it C}$ given by 
$X\tilde{\otimes}Y$ and $Y\tilde{\otimes}X$ are naturally equivalent; namely 
there exist functorial isomorphisms $\Psi_{X,Y}:X\tilde{\otimes}Y\rightarrow$
$Y\tilde{\otimes}X$. From this, the Clifford category is induced.\

A quasitriangular Hopf algebra is a pair $(H,R)$ where $H$ is a Hopf algebra 
and $R\in H\otimes H$ is invertible and obeys ($\Delta\otimes id)R$=
$R_{13}R_{23}$, $(id\otimes \Delta )R$=$R_{13}R_{12}$, $\tau\circ\Delta h$=
$R(\Delta h)R^{-1}$, $\forall h\in H$. Here $\tau:V_1\otimes V_2\rightarrow$
$V_2\otimes V_1$ is the twist map. $(H,R)$ is called triangular if, in 
addition, $\tau (R^{-1})$=$R$. The notation used is $R=\sum R^{(1)}\otimes R^{(2)}$. 
For $({\it H},{\it R})$ as a triangular Hopf algebra 
$\Psi_{V_1,V_2}:V_1\tilde{\otimes}V_2\rightarrow V_2\tilde{\otimes}V_1$ and 
$\Psi (v_1\tilde{\otimes}v_2)$=$\tau\circ 
(\rho_1\otimes\rho_2)({\it R})(v_1\tilde{\otimes}v_2)$ makes 
$_{{\it H}}{\it M}$ into a tensor category which, modulo $I(Q)$, corresponds 
to the Clifford category $\cl_Q$.\

A tensor category 
$({\it C},\tilde{\otimes},\phi,\Psi,\underline{1})$ has an object called 
``internal hom" if the contravariant functors 
$F_{X,Y}=Mor(\tilde{\otimes}X,Y)$(i.e. that send $Z$ to the set 
$Mor(Z\tilde{\otimes}X,Y))$ are each representable. In this case the 
representing object in ${\it C}$, the internal hom, is denoted 
$\underline{Hom}(X,Y)$.\

Let {\it H} be a triangular Hopf algebra with bijective antipode, then its 
finite dimensional algebra representations $({\it M}^{f.d.},\tilde{\otimes},
\underline{1})$ is a rigid tensor category (therefore this induces a Clifford
tensor category) with internal hom defined by $\underline{Hom}(V_1,V_2)$=
$Lin_k(V_1,V_2)$ and the map $\hat{\psi}\in Mor(Z,\underline{Hom}(X,Y))$ for
which $(\hat{\psi}\tilde{\otimes}id)\circ ev_{X,Y}$=$\psi$ where 
$ev_{X,Y}:\underline{Hom}(X,Y)\tilde{\otimes}X\rightarrow Y$ is the evaluation 
map. $\underline{Hom}(X,Y)$ is like ``linear maps from $X$ to $Y$" and 
$ev_{X,Y}$ ``applies" this to an element of $X$ to obtain an element of 
$Y$.\

Suppose that ${\it C}$ is a small abelian ($k$-linear) rigid tensor category, and $F$ is a ($k$-linear, exact, faithful) monoidal functor. Then, 
essentialy, there is a triangular $H$ such that ${\it C}$ is equivalent to
$_{{\it H}}{\it M}$. If in addition $F$ coincides with the 
commutativity constrain in $Vec$. Then the reconstructed $H$ is cocommutative 
i.e. essentially of the form $kG$. Here $G$ is a group and $kG$ denote the
vector space with basis $G$.\

We can make the axioms of a rigid quasitensor category even more explicit by
choosing basis for each of the spaces $V_j$. Thus, let $\{e^j_m\}$ be a basis
for $V_j$, then the general decomposition $V_{j_1}\tilde{\otimes}V_{j_2}\cong$
$\tilde{\oplus}_j\gamma^j_{j_1j_2}\tilde{\otimes}V_j$ takes the form
$$
e^{j_1}_{m_1}\tilde{\otimes}e^{j_2}_{m_2}=
\sum_{j,m}
\left[
\begin{array}{ccc}
j_1   & j_2  & j\\
m_1   & m_2  & m
\end{array}
\right]
\tilde{\otimes}\;\;
e^{j(j_1j_2)}_m
$$
where 
$
\left[
\begin{array}{ccc}
j_1   & j_2  & j\\
m_1   & m_2  & m
\end{array}
\right]\in \gamma^j_{j_1j_2}
$
are called the generalized Clebsch-Gordan coefficients(CGC).\

Recall that in the ordinary Lorentz group, the spinor and vector 
representations are connected by the corresponding Clebsch-Gordan
coefficients, or $\gamma$ matrix in a more familiar terminology. The superfix
$^{(j_1j_2)}$ on the vectors $e^j_m$ is to remind us that they are being viewed 
in $V_{j_1}\tilde{\otimes}V_{j_2}$ according to the isomorphism above 
described. Thus, CGC transform the basis 
$\left\{e^{j(j_1j_2)}_m\right\}$ to the standard basis 
$\left\{e^{j_1}\tilde{\otimes}e^{j_2}\right\}$ of 
$V_{j_1}\tilde{\otimes}V_{j_2}$.\

We can take advantage of this and define $q$-deformed Fierz identities also. 
Let $\psi_a$ be a Majorana $q$-spinor \cite{Ca}; namely 
$\psi_a$=
$\left(
\begin{array}c
Z^k_a\\
(\bar{Z}_a\varepsilon^{-1})^{\bar{k}}
\end{array}
\right)$,
where $a=1,2$ and $Z_a$ is a q-spinor and $\varepsilon$ a metric \cite{Ca}. Then 
we define the following currents;
$$
J_q=\frac{1}{q\sqrt{Q}}\bar{\psi}_1\psi_2\;\;,\;\;\;
J^{\mu}_q=\frac{1}{q\sqrt{Q}}\bar{\psi}_1\gamma^{\mu}_q\psi_2\;\;,\;\;\;
J^{5}_q=\frac{1}{q\sqrt{Q}}\bar{\psi}_1\gamma^{5}_q\psi_2,
$$
$$
J^{\mu\nu}_q=\frac{1}{q\sqrt{Q}}\bar{\psi}_1\gamma^{\mu}_q\gamma^{\nu}_q\psi_2
\;\;\mbox{  and  }\;\;\;
J^{\mu\nu\tau}_q=\frac{1}{q\sqrt{Q}}\bar{\psi}_1\gamma^{\mu}_q\gamma^{\nu}_q
\gamma^{\tau}_q\psi_2,
$$
where $\mu, \nu, \tau=0,+,-,3$ correspond to the indices in $\gamma^{\mu}_q$
previously defined, and $\gamma^5_q=
\gamma^0_q\gamma^+_q\gamma^-_q\gamma^3_q$.\

There exists a finite number of well defined relationships among these 
$q$-Fierz identities. We can easily find the following relations among 
them;
$$
J^{53}_q=-q^2J^{50}_q\;\;,\;\;\;J^{0-}_q=-J^{+3}_q\;\;,\;\;\;
J^{35}_q=J^{05}_q\;\;,\;\;\;J^{-0}_q=q^{-2}J^{3+}_q,
$$
$$
J^{0+}_q=q^2J^{-3}_q\;\;,\;\;\;J^{5+}_q=J^{+-}_q\;\;\mbox{ and }\;\;\;
J^{+0}_q=-J^{3-}_q.
$$
Even more, some other relations appear among these $q$-deformed currents. For
 example $q^4J^2_q-(J^{03}_q)^2$=$Q(1-q^{-4})(J^{5}_q)^2$ where the following particular 
commutation relation between the $q$-spinor and its adjoint has been imposed;
$Z(\epsilon {\bar Z})$=$k{\hat R}_q(\epsilon{\bar Z})Z$, being $k$ a real
number and ${\hat R}_q$=$q\sum_{\rho}e^{\rho}_{\rho}\otimes e^{\rho}_{\rho}+
\newline
\sum_{\rho\neq\sigma}e^{\sigma}_{\rho}\otimes e^{\rho}_{\sigma}+(q-q^{-1})
\sum_{\rho < \sigma}e^{\rho}_{\rho}\otimes e^{\sigma}_{\sigma}$, where
$\rho,\sigma=1,2$ and $e^{\sigma}_{\rho}$ is the basis of the matrix. This 
comes as a natural consequence of the particular braided algebra (reflection
equation) chosen \cite{Az}.\ 

\section{Summary and Conclusions}
Quantum algebras are relatively new mathematical structures which provide an 
exciting generalization of the concept of symmetry. Quantum spin chains 
provide the simplest examples of physical systems which have a quantum algebra 
as invariance.\

Assume that $\cl_{3,1}$ is involved in a theory describing a physical 
reality. We know that any physical theory describes well only a limited class
of phenomena, for the phenomena beyond this class one must modify the theory.
In certain cases such a modification consists in introducing one fundamental 
constant $q$ (small parameter) in the new more general theory (f.e. a quantum 
theory). Within this new theory, $\cl_{3,1}$ retains its validity only in 
the approximate sense. The old theory can be recovered in a limit value for 
$q$. Studying all possible deformations of spacetime Clifford algebras one 
may discover ways leading to more general theories that might better describe 
the reality.\

In this context we are presenting a short review on the theory representation of Hopf
algebras on Clifford algebras. Some examples are also given.\

We intend to search ${\it C}_q(3,1)$ as a symmetry acting on Clifford 
spacetime paving the way to a new approach to quantized spacetime. It turns 
out that upon acting on Clifford spacetime, quantum (i.e. $q$-deformed) 
Clifford spacetime algebra induces a change of metric. From this, we assume 
any quantum symmetry acting on spacetime as performed by a map whose
action can be expressed in terms of well defined fluctuations of metric.\

In the second example we propose $\widehat{CH}_q(2)$ acting as a symmetry of 
su(2); i.e. as a quantum symmetry of the isospin space. Again the result can 
be expressed in terms of particular fluctuations in the metric of the Lie 
algebra.\

Actually these two cases remind us of the map which can be exhibited 
between the generators obeying $SU_q(2)$ commutation relations and the ones 
satisfying su(2) algebra. The properly modified open anisotropic Hiesenberg 
chain posseses $SU_q(2)$ invariance as a higher non-manifest symmetry. We are 
thinking of ${\it C}_q(3,1)$ as a higher non-manifest symmetry of the 
spacetime that realizes quantum symmetries of Clifford algebras.\

As done in the multivectorial generalization of the Cartan map, we propose 
comodule maps to explicitly show isotopic spaces where this symmetry is 
realized. Here we propose to study the category of Clifford representations 
of $q$-deformed (in some cases also Hopf algebras) structures.\

In this paper we show how the Fierz identities, written in terms of products
of generalized multivectorial Cartan maps, have ${\it C}_q(3,1)$ as a higher
non-manifest symmetry as well as su(2) has $\widehat{CH}_q(2)$ as a 
non-manifest symmetry also.\

Since the $c<1$ unitary rational conformal theories can be projected out from 
the $SU_q(2)$-invariant spin 1/2 chain for $q$ a primitive root of unity (in
the thermodynamic limit) we expect ${\it C}_q(3,1)$-invariant Fierz identities 
to be an extension of sensitive meaningful magnitudes in
field theory. Even more, the representation theory of $SU_q(2)$ has 
interesting connections to the solutions of closed spin 1/2 anisotropic Heisenberg chain. 
We think that the representation theory of ${\it C}_q(3,1)$ may have 
connections to some interesting lattice model for spacetime Clifford 
algebra. Some $q$-deformed Fierz identities are also studied.\

The spin 1 Faddeev-Zamolodchikov chain has a class of boundary terms, which 
makes the model $SU_q(2)$-symmmetric. We wonder what kind of discrete model 
with a class of boundary terms corresponds to ${\it C}_q(3,1)$-symmetric Fierz 
identities. It would be interesting to find physical systems with this 
symmetry.\

It has been shown \cite{Cal} that a Hamiltonean is actually already invariant 
under the classical group by virtue of its invariance under the quantum group. 
Since Fierz identities are ${\it C}_q(3,1)$ invariant they certainly are 
$\cl_{3,1}$ invariant. This last symmetry can be realized in a variety of 
ways, in particular the results reported by Rodr\'{\i}guez-Romo \cite{Su1}, 
\cite{Su3}, \cite{Su4}.\

Generalized statistics of physical particles are closely connected with the 
invariance under quantum groups. This invariance provides the possibility to 
construct parafermions possesing generalized statistics which interpolates 
the physical particles \cite{Sm}. This means in our case that the Fierz 
identities already constructed for fermions can be generalized to parafermions in
a straightforward manner.\

\section{Acknowledgments}
This research was partially supported by CONACYT-M\'exico, Ref 4336-E.\

\section{Figure Captions}
Figure 1. Axioms for left $H$-module and $H$-comodule structures written as
diagrams.

\newpage

\begin{picture}(350,250)(0,0)
\put(230,160){\vector(0,1){45}}
\put(230,150){V}
\put(265,150){\vector(1,0){40}}
\put(235,180){$\beta$}
\put(310,180){$\Delta\otimes id$}
\put(220,210){${\it H}\otimes V$}
\put(270,220){$id\otimes \beta$}
\put(270,140){$\beta$}
\put(265,210){\vector(1,0){40}}
\put(320,150){${\it H}\otimes V$}
\put(350,160){\vector(0,1){45}}
\put(310,210){${\it H}\otimes{\it H}\otimes V$}

\put(30,205){\vector(0,-1){45}}
\put(70,140){$\alpha$}
\put(70,220){$id\otimes \alpha$} 
\put(35,180){$\alpha$}
\put(115,180){$\cdot\otimes id$}
\put(30,150){V}
\put(105,150){\vector(-1,0){40}}
\put(20,210){${\it H}\otimes V$}
\put(105,215){\vector(-1,0){40}}
\put(120,150){${\it H}\otimes V$}
\put(150,205){\vector(0,-1){45}}
\put(110,215){${\it H}\otimes{\it H}\otimes V$}

\put(220,110){${\it H}\otimes V$} 
\put(230,40){$V$}  
\put(320,100){$V$} 
\put(235,70){$\beta$} 
\put(270,110){$\epsilon\otimes id$}   
\put(230,50){\vector(0,1){45}}  
\put(265,100){\vector(1,0){40}} 
\put(250,50){\vector(1,1){45}}
\put(290,78){$1_{\it H}\otimes id$}
 
\put(30,40){$V$} 
\put(10,110){${\it H}\otimes V$}     
\put(120,105){$V$}  
\put(70,110){$\eta\otimes id$} 
\put(30,95){\vector(0,-1){45}} 
\put(105,105){\vector(-1,0){40}}    
\put(35,70){$\alpha$}  
\put(105,105){\vector(-1,-1){45}}          
\put(110,70){$1_{\it H}\otimes id$}
\end{picture}\newline
\hspace*{6cm}Figure 1.


\newpage
\begin{thebibliography}{93}
\bibitem{Dr} V.G. Drinfeld, {\it Sov. Math. Dokl.} 27 (1983), 68.\\
V.G. Drinfeld, {\it Sov. Math. Dokl.} 28 (1983), 667.\\
V.G. Drinfeld, {\it Sov. Math. Dokl.} 32 (1985), 254.\\
V.G. Drinfeld, {\it Funct. Anal. Appl.} 20 (1986), 69.
\bibitem{So} F. Sommen and M. Wotkins, ``Introducing $q$-deformation on the
level of vector variables", to appear.
\bibitem{Ma} Yu. I. Manin, {\it Ann. Inst. Fourier} (Grenoble) 37 (1987), 
191.\\
Yu I. Manin, {\it Quantum groups and non-commutative geometry}. Centre de 
Recherches Math\'ematiques. Universit\'e de Montreal (1988).
\bibitem{Ch} K. Fujikawa, L.C. Kwek and C.H. Oh, UT-702(1995), NUS-HEP-95-05.\\
M. Fitchm\"uller, A. Lorek and J. Wess, MPI-PhT/95-109.\\
S. Rodr\'{\i}guez-Romo, {\it XXth International Conference on Differential
Geometric Methods in Theoretical Physics}, S.Catto (Ed),1 (1992), 603.\\
S. Rodr\'{\i}guez-Romo, {\it Geometric and Quantum Aspects of Integrable
Systems}, G.F. Helminck (Ed), Lecture Notes in Physics 424 (1993), 213.\\ 
Z.Chang, J-X. Wang, H. Yan, {\it J. Math. Phys.} 32 (1991), 3241.\\
J.A. Minahan, {\it Mod. Phys. Lett A}, 5 (1990), 2625.\\
H. Yan, J.Phys. A, 23 (1990), L1155.
\bibitem {Sm} F.A. Smirnov, {\it Commun. Math. Phys.} 132 (1990), 415.
\bibitem{Ca} U. Carow-Watamura, M.Schlieker, M. Scholl, S. Watamura, 
{\it Z. Phys.C }48 (1990), 159.\\
U. Carow-Watamura, M.Schlieker, M. Scholl, S. Watamura, 
{\it Int. Jour. of Mod. Phys. A} 6 (1991), 3081.
\bibitem {Su1} S. Rodr\'{\i}guez-Romo, F. Viniegra and J. Keller, ``Clifford
Algebras and their Applications in Mathematical Physics", A. Micali Ed.,
Kluwer Academic Publisher, Netherlands (1992), 479.\\
Suemi Rodr\'{\i}guez-Romo, {\it Found. of Phys.} 23(1993), 1535.
\bibitem{Su3} Jaime Keller and Suemi Rodr\'{\i}guez-Romo, {\it J. Math.Phys.}
32 (1991), 1591.
\bibitem{Su4} Jaime Keller and Suemi Rodr\'{\i}guez-Romo, {\it J. Math. Phys.}
31 (1990), 2501.
\bibitem{Cu} R. Cuerno, C. G\'omez, E. L\'opez and G. Sierra, 
{\it Phys. Lett.B} 307 (1993), 56.
\bibitem{Re1} F. Reifler, {\it J. Math. Phys.} 25 (1984), 1088.\\
F. Reifler, {\it J. Math. Phys.} 26 (1985), 542.\\
F. Reifler, {\it J. Math. Phys.} 26 (1985), 2059.\\
F. Reifler, {\it J. Math. Phys.} 27 (1986), 2803.
\bibitem{Az} J.A. de Azc\'arraga and F. Rodenas, q-alg-9510011.
\bibitem{Cal} D.G. Caldi,`` Quantum Groups", Proccedings of Argonne Workshop,
T. Curtright, David Fairlie and C. Zachos Eds, Word Scientific,  1991.
\end{thebibliography}
\end{document}